\documentclass{article}
\usepackage{amsfonts}

\input{tcilatex}
\begin{document}

\begin{center}
\textbf{DISCRETE FRACTIONAL STURM-LIOUVILLE EQUATIONS}

\bigskip

Erdal BAS$^{a}$ and Ramazan OZARSLAN$^{a\ast }$

\bigskip

$^{a}$\textit{Firat University, Science Faculty, Department of Mathematics,
23119 Elazig/Turkey}

\bigskip

\textit{e-mail: erdalmat@yahoo.com, }$^{\ast }$\textit{%
ozarslanramazan@gmail.com}

\bigskip
\end{center}

\noindent \textbf{Abstract}. In this study, we define discrete fractional
Sturm-Liouville (DFSL) operators within Riemann-Liouville and Gr\"{u}%
nwald-Letnikov fractional operators with both delta and nabla operators. We
show self-adjointness of the DFSL operator for the first time and prove some
spectral properties, like orthogonality of distinct eigenfunctions, reality
of eigenvalues, paralelly in integer and fractional order differential
operator counterparts.

\bigskip


\bigskip


\bigskip

\begin{center}
\textbf{1. Introduction}

\bigskip
\end{center}

Fractional calculus is rather attractive subject due to having wide-ranging
application areas of theoretical and applied sciences. Although there is a
large number of a worthwhile mathematical works on the fractional
differential calculus, there is no noteworthy parallel improvement of
fractional difference calculus up to lately. This improvement has shown that
discrete fractional calculus has certain unforeseen hardship.

Fractional sums and differences were obtained firstly in Diaz-Osler \cite%
{diaz}, Miller and Ross \cite{miller} and Gray and Zhang \cite{gray} and
they found discrete types of fractional integrals and derivatives. Later,
several authors started to touch upon discrete fractional calculus \cite%
{goodrich, dumitru, lizama, mohan, caputo}. Nevertheless, discrete
fractional calculus is a rather novel area. The first studies has been done
by Atici et al. $\left[ 2-5\right] $, Anastassiou \cite{anas, anas2},
Abdeljawad et al. $\left[ 9-13\right] $, and Cheng \cite{jin}, and so forth.

Self-adjoint operators have an important place in differential operators.
Especially, various spectral properties of Sturm-Liouville differential
operators, like orthogonality of distinct eigenfunctions, reality of
eigenvalues, could be analyzed by the help of the self-adjointness.
Self-adjointness of difference operators enabled to be analyzed integer
order Sturm-Liouville difference operators parallelly to differential
counterpart.

Self-adjointness of fractional Sturm-Liouville differential operators have
been shown by $\left[ 23-26\right] $. However, although self-adjointness of
DFSL operators has been mentioned, it has not been proved yet. In this
study, we define DFSL operators, which is defined differently from \cite%
{almeida}, and prove the self-adjointness for the first time. Hereupon, we
prove some spectral properties of the operator paralelly in integer and
fractional order differential counterparts. More recently, Almeida et al. 
\cite{almeida} proved the orthogonality of distinct eigenfunctions, reality
of eigenvalues of DFSL operator and we prove these properties by using
self-adjointness.

In this study, differently we analyze DFSL equation within Riemann-Liouville
and Gr\"{u}nwald-Letnikov fractional operators with both delta and nabla
operators. The aim of this paper is to contribute to the theory of DFSL
operator.

In this paper, we discuss DFSL equation in two different ways;

\noindent $i)$ (nabla right and left) Riemann-Liouville (\textbf{R-L)}
fractional operator,%
\[
L_{1}x\left( t\right) =\nabla _{a}^{\mu }\left( p\left( t\right) _{b}\nabla
^{\mu }x\left( t\right) \right) +q\left( t\right) x\left( t\right) =\lambda
r\left( t\right) x\left( t\right) , 
\]

\noindent $ii)$ (delta left and right) Gr\"{u}nwald-Letnikov (\textbf{G-L)}
fractional operator,%
\[
L_{2}x\left( t\right) =\Delta _{-}^{\mu }\left( p\left( t\right) \Delta
_{+}^{\mu }x\left( t\right) \right) +q\left( t\right) x\left( t\right)
=\lambda r\left( t\right) x\left( t\right) . 
\]

\begin{center}
\bigskip

\textbf{2. Preliminaries}

\bigskip
\end{center}

\noindent \textbf{Definition 1.} \cite{goodrich} Delta and nabla difference
operators are defined by respectively%
\begin{equation}
\Delta x\left( t\right) =x\left( t+1\right) -x\left( t\right) ,\qquad \nabla
x\left( t\right) =x\left( t\right) -x\left( t-1\right) .  \tag{1}  \label{1}
\end{equation}

\bigskip

\noindent \textbf{Definition 2. }\cite{bohner} Falling function is defined
by, $\alpha \in 
%TCIMACRO{\U{211d} }%
%BeginExpansion
\mathbb{R}
%EndExpansion
,$ 
\begin{equation}
t^{\underline{\alpha }}=\frac{\Gamma \left( \alpha +1\right) }{\Gamma \left(
\alpha +1-n\right) },  \tag{2}  \label{2}
\end{equation}%
where $\Gamma $ is the gamma function.

\bigskip

\noindent \textbf{Definition 3. }\cite{bohner} Rising function is defined
by, $\alpha \in 
%TCIMACRO{\U{211d} }%
%BeginExpansion
\mathbb{R}
%EndExpansion
,$ 
\begin{equation}
t^{\bar{\alpha}}=\frac{\Gamma \left( t+\alpha \right) }{\Gamma \left(
t\right) }.  \tag{3}  \label{3}
\end{equation}%
Gamma function must be well-defined in the above two definitions.

\bigskip

\noindent \textbf{Remark 1. }Delta and nabla operators have the following
properties%
\begin{eqnarray}
\Delta t^{\underline{\alpha }} &=&\alpha t^{\underline{\alpha -1}}, 
\TCItag{4}  \label{4} \\
\nabla t^{\bar{\alpha}} &=&\alpha t^{\overline{\alpha -1}}.  \nonumber
\end{eqnarray}

\bigskip

\noindent \textbf{Definition 4. }\cite{ferhan, thabet, miller} Fractional
sum operators are defined by,

\noindent $\left( i\right) $ The nabla left fractional sum of order $\mu >0$
is defined 
\begin{equation}
\nabla _{a}^{-\mu }x\left( t\right) =\frac{1}{\Gamma \left( \mu \right) }%
\sum_{s=a+1}^{t}\left( t-\rho \left( s\right) \right) ^{\overline{\mu -1}%
}x\left( s\right) ,\text{ }t\in 
%TCIMACRO{\U{2115} }%
%BeginExpansion
\mathbb{N}
%EndExpansion
_{a+1},  \tag{5}  \label{6}
\end{equation}%
\noindent $\left( ii\right) $ The nabla right fractional sum of order $\mu
>0 $ is defined 
\begin{equation}
_{b}\nabla ^{-\mu }x\left( t\right) =\frac{1}{\Gamma \left( \mu \right) }%
\sum_{s=t}^{b-1}\left( s-\rho \left( t\right) \right) ^{\overline{\mu -1}%
}x\left( s\right) ,\text{ }t\in \text{ }_{b-1}%
%TCIMACRO{\U{2115} }%
%BeginExpansion
\mathbb{N}
%EndExpansion
,  \tag{6}  \label{7}
\end{equation}%
where $\rho \left( t\right) =t-1$ is called backward jump operators, $%
%TCIMACRO{\U{2115} }%
%BeginExpansion
\mathbb{N}
%EndExpansion
_{a}=\left \{ a,a+1,...\right \} ,$ $_{b}%
%TCIMACRO{\U{2115} }%
%BeginExpansion
\mathbb{N}
%EndExpansion
=\left \{ b,b-1,...\right \} $.

\bigskip

\noindent \textbf{Definition 5. }\cite{thabet2, thabet4} Fractional
difference operators are defined by,

\noindent $\left( i\right) $ The nabla left fractional difference of order $%
\mu >0$ is defined 
\begin{equation}
\nabla _{a}^{\mu }x\left( t\right) =\nabla _{a}^{n}\nabla _{a}^{-\left(
n-\mu \right) }x\left( t\right) =\frac{\nabla ^{n}}{\Gamma \left( n-\mu
\right) }\sum_{s=a+1}^{t}\left( t-\rho \left( s\right) \right) ^{\overline{%
n-\mu -1}}x\left( s\right) ,\text{ }t\in 
%TCIMACRO{\U{2115} }%
%BeginExpansion
\mathbb{N}
%EndExpansion
_{a+1},  \tag{7}  \label{9}
\end{equation}%
$\noindent \left( ii\right) $ The nabla right fractional difference of order 
$\mu >0$ is defined 
\begin{equation}
_{b}\nabla ^{\mu }x\left( t\right) =\left( -1\right) ^{n}\nabla
_{b}^{n}\nabla _{a}^{-\left( n-\mu \right) }x\left( t\right) =\frac{\left(
-1\right) ^{n}\Delta ^{n}}{\Gamma \left( n-\mu \right) }\sum_{s=a+1}^{t}%
\left( s-\rho \left( t\right) \right) ^{\overline{n-\mu -1}}x\left( s\right)
,\text{ }t\in \text{ }_{b-1}%
%TCIMACRO{\U{2115} }%
%BeginExpansion
\mathbb{N}
%EndExpansion
.  \tag{8}  \label{10}
\end{equation}%
Fractional differences in $\left( \ref{9}-\ref{10}\right) $ are called the 
\textbf{Riemann--Liouville (R-L) }definition of the $\mu $-th order nabla
fractional difference.

\bigskip

\noindent \textbf{Definition 6. }$\left[ 18-20\right] $ Fractional
difference operators are defined by,

\noindent $\left( i\right) $ The delta left fractional difference of order $%
\mu ,$ $0<\mu \leq 1,$ is defined%
\begin{equation}
\Delta _{-}^{\mu }x\left( t\right) =\frac{1}{h^{\mu }}\sum_{s=0}^{t}\left(
-1\right) ^{s}\frac{\mu \left( \mu -1\right) ...\left( \mu -s+1\right) }{s!}%
x\left( t-s\right) ,\text{ }t=1,...,N.  \tag{9}  \label{11}
\end{equation}%
$\noindent \left( ii\right) $ The delta right fractional difference of order 
$\mu ,$ $0<\mu \leq 1,$ is defined%
\begin{equation}
\Delta _{+}^{\mu }x\left( t\right) =\frac{1}{h^{\mu }}\sum_{s=0}^{N-t}\left(
-1\right) ^{s}\frac{\mu \left( \mu -1\right) ...\left( \mu -s+1\right) }{s!}%
x\left( t+s\right) ,\text{ }t=0,..,N-1,  \tag{10}  \label{12}
\end{equation}%
Fractional differences in $\left( \ref{11}-\ref{12}\right) $ are called the 
\textbf{Gr\"{u}nwald--Letnikov (G-L) }definition of the $\mu $-th order
delta fractional difference.

\bigskip

\noindent \textbf{Definition 7. }\cite{thabet4} We define the integration by
parts formula for\textbf{\ R-L} nabla fractional difference operator, $u$ is
defined on $_{b}%
%TCIMACRO{\U{2115} }%
%BeginExpansion
\mathbb{N}
%EndExpansion
$ and $v$ is defined on $%
%TCIMACRO{\U{2115} }%
%BeginExpansion
\mathbb{N}
%EndExpansion
_{a}$, then 
\begin{equation}
\sum_{s=a+1}^{b-1}u\left( s\right) \nabla _{a}^{\mu }v\left( s\right)
=\sum_{s=a+1}^{b-1}v\left( s\right) _{b}\nabla ^{\mu }u\left( s\right) . 
\tag{11}  \label{16}
\end{equation}

\bigskip

\noindent \textbf{Definition 8. }\cite{almeida, bourdin} We define the
integration by parts formula for \textbf{G-L }delta fractional difference
operator, $u$, $v$ is defined on $\left \{ 0,1,...,n\right \} $, then 
\begin{equation}
\sum_{s=0}^{n}u\left( s\right) \Delta _{-}^{\mu }v\left( s\right)
=\sum_{s=0}^{n}v\left( s\right) \Delta _{+}^{\mu }u\left( s\right) . 
\tag{12}  \label{17}
\end{equation}

\begin{center}
\bigskip

\textbf{3. Main Results}
\end{center}

\bigskip

\noindent \textbf{3.1. Discrete Fractional Sturm-Liouville Equations}

\bigskip

We consider DFSL equations in two different ways;

\noindent $i)$ (nabla right and left) \textbf{R-L} fractional operator is
defined by,%
\begin{equation}
L_{1}x\left( t\right) =\nabla _{a}^{\mu }\left( p\left( t\right) _{b}\nabla
^{\mu }x\left( t\right) \right) +q\left( t\right) x\left( t\right) =\lambda
r\left( t\right) x\left( t\right) ,  \tag{13}  \label{23}
\end{equation}%
where $p\left( t\right) >0,$ $r\left( t\right) >0,$ $q\left( t\right) $ is
defined and real valued, $\lambda $ is the spectral parameter, $x\left(
t\right) \in l^{2}\left[ a+1,b-1\right] .$

\noindent $ii)$ (delta left and right) \textbf{G-L} fractional operator is
defined by,%
\begin{equation}
L_{2}x\left( t\right) =\Delta _{-}^{\mu }\left( p\left( t\right) \Delta
_{+}^{\mu }x\left( t\right) \right) +q\left( t\right) x\left( t\right)
=\lambda r\left( t\right) x\left( t\right) ,  \tag{14}  \label{27}
\end{equation}%
where $p,q,r,\lambda $ is as defined above, $x\left( t\right) \in l^{2}\left[
0,n\right] .$

\bigskip

Firstly, let's consider the equation $\left( \ref{23}\right) $ and give the
following theorems and proofs$,$

\bigskip

\noindent \textbf{Theorem 9. }DFSL operator $L_{1}$, denoted by the equation 
$\left( \ref{23}\right) ,$ is self-adjoint.

\bigskip

\textbf{Proof.}%
\begin{eqnarray}
u\left( t\right) L_{1}v\left( t\right) &\mathbf{=}&u\left( t\right) \nabla
_{a}^{\mu }\left( p\left( t\right) _{b}\nabla ^{\mu }v\left( t\right)
\right) +u\left( t\right) q\left( t\right) v\left( t\right) ,  \TCItag{15} \\
v\left( t\right) L_{1}u\left( t\right) &\mathbf{=}&v\left( t\right) \nabla
_{a}^{\mu }\left( p\left( t\right) _{b}\nabla ^{\mu }u\left( t\right)
\right) +v\left( t\right) q\left( t\right) u\left( t\right) .  \TCItag{16}
\end{eqnarray}%
If $\left( 15-16\right) $ is subtracted from each other%
\[
u\left( t\right) L_{1}v\left( t\right) -v\left( t\right) L_{1}u\left(
t\right) =u\left( t\right) \nabla _{a}^{\mu }\left( p\left( t\right)
_{b}\nabla ^{\mu }v\left( t\right) \right) -v\left( t\right) \nabla
_{a}^{\mu }\left( p\left( t\right) _{b}\nabla ^{\mu }u\left( t\right)
\right) 
\]%
and definite sum operator to the both side of the last equality is applied,
we have%
\begin{equation}
\sum_{s=a+1}^{b-1}\left( u\left( s\right) L_{1}v\left( s\right) -v\left(
s\right) L_{1}u\left( s\right) \right) =\sum_{s=a+1}^{b-1}u\left( s\right)
\nabla _{a}^{\mu }\left( p\left( s\right) _{b}\nabla ^{\mu }v\left( s\right)
\right) -\sum_{s=a+1}^{b-1}v\left( s\right) \nabla _{a}^{\mu }\left( p\left(
s\right) _{b}\nabla ^{\mu }u\left( s\right) \right) .  \tag{17}
\end{equation}%
If we apply the integration by parts formula in $\left( \ref{16}\right) $ to
right hand side of $\left( 17\right) ,$ we have%
\begin{eqnarray*}
\sum_{s=a+1}^{b-1}\left( u\left( s\right) L_{1}v\left( s\right) -v\left(
s\right) L_{1}u\left( s\right) \right) &=&\sum_{s=a+1}^{b-1}p\left( s\right)
_{b}\nabla ^{\mu }v\left( s\right) _{b}\nabla ^{\mu }u\left( s\right) \\
&&-\sum_{s=a+1}^{b-1}p\left( s\right) _{b}\nabla ^{\mu }u\left( s\right)
_{b}\nabla ^{\mu }v\left( s\right) \\
&=&0,
\end{eqnarray*}%
\[
\left \langle L_{1}u,v\right \rangle =\left \langle u,L_{1}v\right \rangle . 
\]%
The proof completes.

\bigskip

\noindent \textbf{Theorem 10. }Eigenfunctions, corresponding to distinct
eigenvalues, of the equation $\left( \ref{27}\right) $ are orthogonal.

\bigskip

\textbf{Proof. }Let $\lambda _{\alpha }$ and $\lambda _{\beta }$ are two
different eigenvalues corresponds to eigenfunctions $u\left( n\right) $ and $%
v\left( n\right) $ respectively for the the equation $\left( \ref{23}\right) 
$,%
\begin{eqnarray*}
\nabla _{a}^{\mu }\left( p\left( t\right) _{b}\nabla ^{\mu }u\left( t\right)
\right) +q\left( t\right) u\left( t\right) -\lambda _{\alpha }r\left(
t\right) u\left( t\right) &=&0, \\
\nabla _{a}^{\mu }\left( p\left( t\right) _{b}\nabla ^{\mu }v\left( t\right)
\right) +q\left( t\right) v\left( t\right) -\lambda _{\beta }r\left(
t\right) v\left( t\right) &=&0,
\end{eqnarray*}%
If we multiply last two equations to $v\left( n\right) $ and $u\left(
n\right) $ respectively, subtract from each other and apply definite sum
operator, since the self-adjointness of the operator $L_{1},$ we get 
\[
\left( \lambda _{\alpha }-\lambda _{\beta }\right) \sum_{s=a+1}^{b-1}r\left(
s\right) u\left( s\right) v\left( s\right) =0, 
\]%
since $\lambda _{\alpha }\neq \lambda _{\beta ,}$%
\begin{eqnarray*}
\sum_{s=a+1}^{b-1}r\left( s\right) u\left( s\right) v\left( s\right) &=&0 \\
\left \langle u\left( t\right) ,v\left( t\right) \right \rangle &=&0
\end{eqnarray*}%
The proof completes..

\bigskip

\noindent \textbf{Theorem 11. }All eigenvalues of the equation $\left( \ref%
{23}\right) $ are real.

\bigskip

\textbf{Proof. }Let $\lambda =\alpha +i\beta ,$ since the self-adjointness
of the operator $L_{1},$ we have%
\begin{eqnarray*}
\left \langle L_{1}u,u\right \rangle &=&\left \langle u,L_{1}u\right \rangle
, \\
\left \langle \lambda ru,u\right \rangle &=&\left \langle u,\lambda ru\right
\rangle ,
\end{eqnarray*}%
\[
\left( \lambda -\overline{\lambda }\right) \left \langle u,u\right \rangle
_{r}=0 
\]%
Since $\left \langle u,u\right \rangle _{r}\neq 0,$ 
\[
\lambda =\overline{\lambda } 
\]%
and hence $\beta =0.$ The proof completes.

\bigskip

Secondly, let's consider the equation $\left( \ref{23}\right) $ and give the
following theorems and proofs$,$

\bigskip

\noindent \textbf{Theorem 12. }DFSL operator $L_{2}$, denoted by the
equation $\left( \ref{27}\right) ,$ is self-adjoint.

\bigskip

\textbf{Proof. }%
\begin{eqnarray}
u\left( t\right) L_{2}v\left( t\right) &\mathbf{=}&u\left( t\right) \Delta
_{-}^{\mu }\left( p\left( t\right) \Delta _{+}^{\mu }v\left( t\right)
\right) +u\left( t\right) q\left( t\right) v\left( t\right) ,  \TCItag{18} \\
v\left( t\right) L_{2}u\left( t\right) &\mathbf{=}&v\left( t\right) \Delta
_{-}^{\mu }\left( p\left( t\right) \Delta _{+}^{\mu }u\left( t\right)
\right) +v\left( t\right) q\left( t\right) u\left( t\right) .  \TCItag{19}
\end{eqnarray}%
If $\left( 18-19\right) $ is subtracted from each other%
\[
u\left( t\right) L_{2}v\left( t\right) -v\left( t\right) L_{2}u\left(
t\right) =u\left( t\right) \Delta _{-}^{\mu }\left( p\left( t\right) \Delta
_{+}^{\mu }v\left( t\right) \right) -v\left( t\right) \Delta _{-}^{\mu
}\left( p\left( t\right) \Delta _{+}^{\mu }u\left( t\right) \right) 
\]%
and definite sum operator to the both side of the last equality is applied,
we have%
\begin{equation}
\sum_{s=0}^{n}\left( u\left( s\right) L_{1}v\left( s\right) -v\left(
s\right) L_{2}u\left( s\right) \right) =\sum_{s=0}^{n}u\left( s\right)
\Delta _{-}^{\mu }\left( p\left( s\right) \Delta _{+}^{\mu }v\left( s\right)
\right) -\sum_{s=0}^{n}v\left( s\right) \Delta _{-}^{\mu }\left( p\left(
s\right) \Delta _{+}^{\mu }u\left( s\right) \right) .  \tag{20}
\end{equation}%
If we apply the integration by parts formula in $\left( \ref{17}\right) $ to
right hand side of $\left( 20\right) ,$ we have%
\begin{eqnarray*}
\sum_{s=0}^{n}\left( u\left( s\right) L_{2}v\left( s\right) -v\left(
s\right) L_{2}u\left( s\right) \right) &=&\sum_{s=0}^{n}p\left( s\right)
\Delta _{+}^{\mu }v\left( s\right) \Delta _{+}^{\mu }u\left( s\right) \\
&&-\sum_{s=0}^{n}p\left( s\right) \Delta _{+}^{\mu }u\left( s\right) \Delta
_{+}^{\mu }v\left( s\right) \\
&=&0
\end{eqnarray*}%
\[
\left \langle L_{2}u,v\right \rangle =\left \langle u,L_{2}v\right \rangle . 
\]%
The proof completes.

\bigskip

\noindent \textbf{Theorem 13. }Eigenfunctions, corresponding to distinct
eigenvalues, of the equation $\left( \ref{27}\right) $ are orthogonal.

\bigskip

\textbf{Proof. }Let $\lambda _{\alpha }$ and $\lambda _{\beta }$ are two
different eigenvalues corresponds to eigenfunctions $u\left( n\right) $ and $%
v\left( n\right) $ respectively for the the equation $\left( \ref{27}\right) 
$,%
\begin{eqnarray*}
\Delta _{-}^{\mu }\left( p\left( t\right) \Delta _{+}^{\mu }u\left( t\right)
\right) +q\left( t\right) u\left( t\right) -\lambda _{\alpha }r\left(
t\right) u\left( t\right) &=&0, \\
\Delta _{-}^{\mu }\left( p\left( t\right) \Delta _{+}^{\mu }v\left( t\right)
\right) +q\left( t\right) v\left( t\right) -\lambda _{\beta }r\left(
t\right) v\left( t\right) &=&0.
\end{eqnarray*}%
If we multiply last two equations to $v\left( n\right) $ and $u\left(
n\right) $ respectively, subtract from each other and apply definite sum
operator, since the self-adjointness of the operator $L_{2},$ we get 
\[
\left( \lambda _{\alpha }-\lambda _{\beta }\right) \sum_{s=0}^{n}r\left(
s\right) u\left( s\right) v\left( s\right) =0, 
\]%
since $\lambda _{\alpha }\neq \lambda _{\beta ,}$%
\begin{eqnarray*}
\sum_{s=0}^{n}r\left( s\right) u\left( s\right) v\left( s\right) &=&0 \\
\left \langle u\left( t\right) ,v\left( t\right) \right \rangle &=&0.
\end{eqnarray*}%
So, the eigenfunctions are orthogonal. The proof completes.

\bigskip

\noindent \textbf{Theorem 14. }All eigenvalues of the equation $\left( \ref%
{27}\right) $ are real.

\bigskip

\textbf{Proof }Let $\lambda =\alpha +i\beta ,$ since the self-adjointness of
the operator $L_{2}$%
\begin{eqnarray*}
\left \langle L_{2}u,u\right \rangle &=&\left \langle u,L_{2}u\right \rangle
, \\
\left \langle \lambda ru,u\right \rangle &=&\left \langle u,\lambda ru\right
\rangle ,
\end{eqnarray*}%
\[
\left( \lambda -\overline{\lambda }\right) \left \langle u,u\right \rangle
_{r}=0 
\]%
Since $\left \langle u,u\right \rangle _{r}\neq 0,$ 
\[
\lambda =\overline{\lambda }, 
\]%
and hence $\beta =0.$ The proof completes.

\bigskip

\begin{center}
\textbf{4. Conclusion}
\end{center}

In this study, differently we analyze DFSL equation within Riemann-Liouville
and Gr\"{u}nwald-Letnikov fractional operators with both delta and nabla
operators. We define DFSL operators, which is defined differently from \cite%
{almeida}, and prove the self-adjointness for the first time. Hereupon, we
prove some spectral properties of the operator paralelly in integer and
fractional order differential counterparts, like orthogonality of distinct
eigenfunctions, reality of eigenvalues. The aim of this paper is to
contribute to the theory of Sturm-Liouville fractional difference operator.

\end{document}